\documentclass[12pt,draft]{amsart}

\usepackage[cp1251]{inputenc}
\usepackage[T2A]{fontenc}
\usepackage[russian, english]{babel}
\usepackage{amsmath,amsfonts,amssymb}
\usepackage{geometry}
\newtheorem{theorem}{Theorem}
\newtheorem{corollary}{Corollary}
\newtheorem{lemma}{Lemma}
\newtheorem{proposition}{Proposition}

\theoremstyle{remark}
\newtheorem{remark}{Remark}

\theoremstyle{definition}

\title{Extension of some Lions-Magenes theorems}

\author{Aleksandr A. Murach}

\address{Institute of Mathematics NAS of Ukraine,
Tereshchenkivska str. 3, Kyiv, Ukraine, 01601; \indent Chernigiv State
Technological University, Shevchenka str. 95, Chernigiv 14027, Ukraine}

\email{murach@imath.kiev.ua}

\subjclass[2000]{Primary 35J40; Secondary 46E35}

\keywords{Elliptic boundary-value problem, generalized solution, Sobolev
spaces, Lions-Magenes theorems, multipliers.}

\dedicatory{Dedicated to the memory of A. Y. Povzner}

\begin{document}

\maketitle

\begin{abstract}
A general form of the  Lions-Magenes theorems on solvability of an elliptic
boundary-value problem in the spaces of nonregular distributions is proved. We
find a general condition on the space of right-hand sides of the elliptic
equation under which the operator of the problem is bounded and has a finite
index on the corresponding couple of Hilbert spaces. Extensive classes of the
spaces satisfying this condition are constructed. They contain the spaces used
by Lions and Magenes.
\end{abstract}

\section{Introduction and statement of the problem}\label{sec1}

Let $\Omega$ be a bounded domain in the Euclidean space $\mathbb{R}^n$,
$n\geq2$, with the boundary $\Gamma$ which is an infinitely smooth closed
manifold of the dimension $n-1$. The domain $\Omega$ is situated locally on the
same side from $\Gamma$.

We consider the nonhomogeneous boundary-value problem in the domain $\Omega$:
\begin{equation}\label{1}
A\,u=f\;\;\text{in}\;\;\Omega,\quad
B_{j}\,u=g_{j}\;\;\text{on}\;\;\Gamma\;\;\text{for}\;\; j=1,\ldots,q.
\end{equation}
In what follows $A$ is a linear differential expression on $\overline{\Omega}$
of an arbitrary even order $2q\geq\nobreak2$, whereas $B_{j}$ with
$j=1,\ldots,q$ is a boundary linear differential expression on $\Gamma$ of
order $m_{j}\leq\nobreak2q-1$. All coefficients of $A$ and $B_{j}$ are
complex-valued functions infinitely smooth on
$\overline{\Omega}:=\Omega\cup\Gamma$ and on $\Gamma$ respectively.

Everywhere in this paper the boundary-value problem \eqref{1} is assumed to be
\textit{regular elliptic}. This means \cite[Ch.~2, Sec.~5.1]{LionsMagenes68}
that the expression $A$ is properly elliptic on $\overline{\Omega}$ and the
collection of boundary expressions $B:=(B_{1},\ldots,B_{q})$ is normal and
satisfies the complementing condition with respect to $A$ on $\Gamma$. It
follows from the condition of normality that all orders $m_{j}$ with
$j=1,\ldots,q$ are mutually distinct.

Along with \eqref{1} we consider the boundary-value problem
\begin{equation}\label{a}
A^{+}\,v=\omega\;\;\text{in}\;\;\Omega,\quad
B^{+}_{j}\,v=h_{j}\;\;\text{on}\;\;\Gamma,\;\;j=1,\ldots,q,
\end{equation}
formally adjoint to the problem \eqref{1} with respect to the Green formula
$$
(Au,v)_{\Omega}+\sum_{j=1}^{q}\;(B_{j}u,\,C_{j}^{+}v)_{\Gamma}
=(u,A^{+}v)_{\Omega}+\sum_{j=1}^{q}\;(C_{j}u,\,B_{j}^{+}v)_{\Gamma}, \quad
u,v\in C^{\infty}(\,\overline{\Omega}\,).
$$
Here, $A^{+}$ is the linear differential expression formally adjoint to $A$ and
having the order $2q$ and coefficients from
$C^{\infty}(\,\overline{\Omega}\,)$. In addition, $\{B^{+}_{j}\}$, $\{C_{j}\}$,
and $\{C^{+}_{j}\}$ are certain normal systems of linear differential boundary
expression with coefficients from $C^{\infty}(\Gamma)$. They orders satisfy the
condition:
$$\mathrm{ord}\,B_{j}+\mathrm{ord}\,C^{+}_{j}=
\mathrm{ord}\,C_{j}+\mathrm{ord}\,B^{+}_{j}=2q-1.
$$
In what follows we denote by $(\cdot,\cdot)_{\Omega}$ and
$(\cdot,\cdot)_{\Gamma}$ the inner products in the spaces $L_{2}(\Omega)$ and
$L_{2}(\Gamma)$ (formed by functions square-integrable over $\Omega$ or
$\Gamma$ respectively) and natural extensions by continuity of these inner
products.

We set
\begin{gather*}
N:=\{u\in C^{\infty}(\,\overline{\Omega}\,):\;Au=0\;\;\mbox{in}\;\; \Omega,\;\;
B_{j}u=0\;\;\mbox{on}\;\;\Gamma\;\;\forall\;\;j=1,\ldots,q\}, \\
N^{+}:=\{v\in C^{\infty}(\,\overline{\Omega}\,):\;A^{+}v=0\;\;\mbox{in}\;\;
\Omega,\;\; B^{+}_{j}v=0\;\;\mbox{on}\;\;\Gamma\;\;\forall\;\;j=1,\ldots,q\}.
\end{gather*}
Since both the problem \eqref{1} and \eqref{a} are regular elliptic, both the
spaces $N$ and $N^{+}$ are finite-dimensional \cite[Ch.~2,
Theorem~5.3]{LionsMagenes68}.

The fundamental property of every elliptic boundary-value problem consists in
that the problem generates the bounded and Fredholm operator on appropriate
couples of functional spaces. Note that a linear bounded operator
$T:E_{1}\rightarrow E_{2}$, with $E_{1}$ and $E_{2}$ being Banach spaces, is
called the Fredholm operator if its kernel $\ker T$ and co-kernel
$\mathrm{coker}\,T:=E_{2}/T(E_{1})$ are both finite-dimensional. The Fredholm
operator $T$ has the closed range in the space $E_{2}$ and the finite index
$\mathrm{ind}\,T:=\ker T-\mathrm{coker}\,T$. This operator naturally generates
the homeomorphism $T:E_{1}/\ker T\leftrightarrow T(E_{1})$.

Let us formulate the classical theorem on elliptic boundary-value problems
(see, e.g. \cite[Ch.~2, Sec.~5.4]{LionsMagenes68}, \cite[Ch.~3,
Sec.~6]{Berezansky68}). In this paper, we restrict ourselves to the Hilbert
spaces case, which is the most important for applications.

\medskip

\noindent\textbf{Theorem 0.} \it The mapping
\begin{equation}\label{b}
u\mapsto(Au,Bu),\quad u\in C^{\infty}(\,\overline{\Omega}\,),
\end{equation}
can be extended by continuity to the bounded and Fredholm operator
\begin{equation}\label{c}
(A,B):\,H^{s+2q}(\Omega)\rightarrow
H^{s}(\Omega)\oplus\bigoplus_{j=1}^{q}\,H^{s+2q-m_{j}-1/2}(\Gamma)
=:\mathcal{H}_{s}(\Omega,\Gamma)
\end{equation}
for every real $s\geq0$. The kernel of this operator coincides with $N$,
whereas the range consists of all vectors
$(f,g_{1},\ldots,g_{q})\in\mathcal{H}_{s}(\Omega,\Gamma)$ satisfying the
condition
\begin{equation}\label{4}
(f,v)_{\Omega}+\sum_{j=1}^{q}\,(g_{j},C^{+}_{j}v)_{\Gamma}=0\quad\mbox{for
every}\quad v\in N^{+}.
\end{equation}
The index of the operator \eqref{c} is equal to $\dim N-\dim N^{+}$ and
independent of $s$. \rm

\medskip

In what follows $H^{\sigma}(\Omega)$ and $H^{\sigma}(\Gamma)$ are Hilbert
spaces with the index $\sigma\in\mathbb{R}$ consisting of some distributions
given in the domain $\Omega$ or on the manifold $\Gamma$ respectively (we will
remind their definitions in Sec. \ref{sec2}). In addition, as usual
$\mathcal{D}'(\Omega)$ and $\mathcal{D}'(\Gamma)$ stand for the linear
topological spaces of all distributions given in $\Omega$ or on $\Gamma$. We
always interpret distributions as antilinear functionals.

Theorem~0 has a \emph{generic} nature because the spaces in which the operator
\eqref{c} acts are common for all elliptic boundary-value problems of the same
order. By this theorem, the operator $(A,B)$ establishes the homeomorphism of
the factor space $H^{s+2q}(\Omega)/N$ onto the subspace
$$
\{(f,g_{1},\ldots,g_{q})\in\mathcal{H}_{s}(\Omega,\Gamma):\,\mbox{\eqref{4} is
true}\}
$$
for each $s\geq0$. Therefore the theorems on operators generated by elliptic
boundary-value problems are called the theorems on homeomorphisms.

Generally, Theorem~0 is not true in the case $s<0$ because the mapping
$u\mapsto B_{j}u$ with $u\in C^{\infty}(\,\overline{\Omega}\,)$ can not be
extended to the bounded operator $B_{j}:H^{s+2q}(\Omega)\rightarrow
\mathcal{D}'(\Gamma)$ if $s+2q\leq m_{j}+1/2$. Therefore we have to use the
space narrower than $H^{s+2q}(\Omega)$ as the domain of $(A,B)$, namely
\begin{equation}\label{2}
D^{s+2q}_{A,X}(\Omega):=\{u\in H^{s+2q}(\Omega):\,Au\in X^{s}(\Omega)\},
\end{equation}
where $X^{s}(\Omega)$ is a Hilbert space imbedded continuously in
$\mathcal{D}'(\Omega)$. In what follows the image $Au$ of
$u\in\mathcal{D}'(\Omega)$ is understood in the theory of distributions. We
endow the space \eqref{2} with the graphics inner product
\begin{equation}\label{d}
(u_{1},u_{2})_{D^{s+2q}_{A,X}(\Omega)}:=(u_{1},u_{2})_{H^{s+2q}(\Omega)}+
(Au_{1},Au_{2})_{X^{s}(\Omega)}
\end{equation}
and the corresponding norm.

The space $D^{s+2q}_{A,X}(\Omega)$ with the inner product \eqref{d} is
complete. Indeed, if $(u_{k})$ is a Cauchy sequence in
$D^{s+2q}_{A,X}(\Omega)$, then by a completeness of $H^{s+2q}(\Omega)$ and
$X^{s}(\Omega)$ there are two limits: $u:=\lim u_{k}$ in
$H^{s+2q}(\Omega)\hookrightarrow\mathcal{D}'(\Omega)$ and $f:=\lim Au_{k}$ in
$X^{s}(\Omega)\hookrightarrow\mathcal{D}'(\Omega)$ (imbeddings are continuous).
Since the differential operator $A$ is continuous in $\mathcal{D}'(\Omega)$, we
deduce from the first limit that $Au=\lim Au_{k}$ in $\mathcal{D}'(\Omega)$.
This implies by the second limit the equality $Au=f\in X^{s}(\Omega)$.
Therefore $u\in D^{s+2q}_{A,X}(\Omega)$ and $\lim u_{k}=u$ in the space
$D^{s+2q}_{A,X}(\Omega)$, that is this space is complete.

J.-L.~Lions and E.~Magenes \cite{LionsMagenes62V, LionsMagenes63VI, Magenes65,
LionsMagenes68} found the certain important examples of $X^{s}(\Omega)$ such
that the mapping \eqref{b} can be extended by continuity to the bounded and
Fredholm operator
\begin{equation}\label{3}
(A,B):\,D^{s+2q}_{A,X}(\Omega)\rightarrow
X^{s}(\Omega)\oplus\bigoplus_{j=1}^{q}\,H^{s+2q-m_{j}-1/2}(\Gamma)
=:\mathcal{X}_{s}(\Omega,\Gamma)
\end{equation}
if $s<0$. In contrast to Theorem~0, the domain of the operator \eqref{3} with
the topology depends on coefficients of the elliptic expression $A$. Therefore
the theorems on properties of the operator \eqref{3} naturally can be termed
the \emph{individual} theorems. Let us formulate two individual theorems proved
by Lions and Magenes.

\medskip

\noindent\textbf{Theorem~LM$_{\mathbf{1}}$} \cite{LionsMagenes62V,
LionsMagenes63VI}. \it Let $s<0$ and $X^{s}(\Omega):=L_{2}(\Omega)$. Then the
mapping \eqref{b} can be extended by continuity to the bounded and Fredholm
operator \eqref{3}. The kernel of this operator coincides with $N$, whereas the
domain consists of all vectors
$(f,g_{1},\ldots,g_{q})\in\mathcal{X}_{s}(\Omega,\Gamma)$ satisfying \eqref{4}.
The index of the operator \eqref{3} is $\dim N-\dim N^{+}$ and independent of
$s$. \rm

\medskip

Here, we especially note the case $s=-2q$, which is important in the spectral
theory of elliptic operators \cite{Grubb68, Grubb96, Mikhailets82,
Mikhailets89}. In this case the space
\begin{equation}\label{e9}
D^{0}_{A,L_{2}}(\Omega)=\{u\in L_{2}(\Omega):\,Au\in L_{2}(\Omega)\}
\end{equation}
is the domain of the maximal operator corresponding to the differential
expression $A$ \cite[Sec.~1.2]{Hermander55}. Even when all coefficients of $A$
are constant, the space \eqref{e9} depends  essentially on each of them. We can
see it from the following result of L.~H\"ormander \cite[Sec.~3.1,
Theorem~3.1]{Hermander55}.

Let both $A_{1}$ and $A_{2}$ be constant-coefficient linear differential
expressions. If $D^{0}_{A_{1},L_{2}}(\Omega)\subseteq
D^{0}_{A_{2},L_{2}}(\Omega)$, then either $A_{2}=\alpha A_{1}+\beta$ for some
$\alpha,\beta\in\mathbb{C}$, or both $A_{1}$ and $A_{2}$ are certain
polynomials in the derivation operator with respect to a vector $e$ and
moreover $\mathrm{ord}\,A_{2}\leq\mathrm{ord}\,A_{1}$. Note that the second
possibility is excluded for elliptic operators.

To formulate the second Lions-Magenes theorem we need the following weighted
space
$$
\varrho H^{s}(\Omega):=\{f=\varrho v:\,v\in H^{s}(\Omega)\,\},
$$
where $s<0$, and a function $\varrho\in C^{\infty}(\Omega)$ is positive. We
endow this space with the inner product
$$
\quad (f_{1},f_{2})_{\varrho H^{s}(\Omega)}:=
(\varrho^{-1}f_{1},\varrho^{-1}f_{2})_{H^{s}(\Omega)}
$$
and the corresponding norm. The space $\varrho H^{s}(\Omega)$ is complete and
imbedded continuously in $\mathcal{D}'(\Omega)$. This follows from that the
operator of multiplication by $\varrho$ is continuous in $\mathcal{D}'(\Omega)$
and establishes the homeomorphism from the complete space $H^{s}(\Omega)$ onto
$\varrho H^{s}(\Omega)$.

We consider a weight function $\varrho:=\varrho_{1}^{-s}$ such that
\begin{equation}\label{f}
\varrho_{1}\in C^{\infty}(\,\overline{\Omega}\,),\quad\varrho_{1}>0\;\;\mbox{in
$\Omega$},\quad\varrho_{1}(x)=\mathrm{dist}(x,\Gamma)\quad\mbox{in a
neighbourhood of $\Gamma$}.
\end{equation}

\medskip

\noindent\textbf{Theorem LM$_{\mathbf{2}}$} \cite[Ch.~2,
Sec.~7.3]{LionsMagenes68}. \it Let $s<0$, $s+1/2\notin\mathbb{Z}$, and
$X^{s}(\Omega):=\rho_{1}^{-s}H^{s}(\Omega)$. Then the assertion of
Theorem~$\mathrm{LM}_{1}$ remains true. \rm

\medskip

We note that Lions and Magenes used a certain Hilbert space $\Xi^{s}(\Omega)$
as $X^{s}(\Omega)$. This space coincides (up to equivalence of norms) with the
weighted space $\rho_{1}^{-s}H^{s}(\Omega)$ for each non half-integer $s<0$
\cite[Ch.~2, Sec.~7.1]{LionsMagenes68}. Theorem~LM$_{2}$ also holds true for
every half-integer $s<0$ if we define the space $X^{s}(\Omega)$ by means of the
complex (holomorphic) interpolation, for instance
\begin{equation}\label{g}
X^{s}(\Omega):=[X^{2s}(\Omega),L_{2}(\Omega)]_{1/2}.
\end{equation}
(See the definition and properties of this interpolation, e.g. in \cite[Ch.~1,
Sec.~14.1]{LionsMagenes68}).

In this paper, we find a general enough condition on the space $X^{s}(\Omega)$
under which the operator \eqref{3} is well defined, bounded, and Fredholm if
$s<0$. The condition consists in the following.

\medskip

\noindent\textbf{Condition I$_{s}$.} The set
$X^{\infty}(\Omega):=X^{s}(\Omega)\cap C^{\infty}(\,\overline{\Omega}\,)$ is
dense in $X^{s}(\Omega)$, and there exists a number $c>0$ such that
\begin{equation}\label{h}
\|\mathcal{O}f\|_{H^{s}(\mathbb{R}^{n})}\leq
c\,\|f\|_{X^{s}(\Omega)}\quad\mbox{for each}\quad f\in X^{\infty}(\Omega).
\end{equation}
Here, $\mathcal{O}f(x):=f(x)$ for $x\in\overline{\Omega}$, and
$\mathcal{O}f(x):=0$ for $x\in\mathbb{R}^{n}\setminus\overline{\Omega}$.

\medskip

In \eqref{h}, we define by $H^{s}(\mathbb{R}^{n})$ the Hilbertian Sobolev space
with index $s$ and given over $\mathbb{R}^{n}$. Note that if $s$ is smaller,
then Condition I$_{s}$ is weaker for the same space $X^{s}(\Omega)$.

Both of the spaces $X^{s}(\Omega):=L_{2}(\Omega)$ and
$X^{s}(\Omega):=\rho_{1}^{-s}H^{s}(\Omega)$ used by Lions and Magenes satisfy
Condition~I$_{s}$.

In this paper, we find all the Hilbertian Sobolev spaces
$X^{s}(\Omega)=H^{\sigma}(\Omega)$ for which Condition~I$_{s}$ is fulfilled. In
addition, we describe the class of all weights $\varrho\in C^{\infty}(\Omega)$
such that the weighted space $X^{s}(\Omega):=\rho H^{s}(\Omega)$ satisfies
Condition~I$_{s}$. This class contains the weight $\rho:=\rho_{1}^{-s}$ as a
particular case. Thus, we get some generalizations of the Lions-Magenes
theorems mentioned above to more extensive classes of the Hilbertian spaces
$X^{s}(\Omega)$ of right-hand sides of the elliptic equation.

Note that we generalize the Lions-Magenes theorems staying in classis of
distributions given in the domain $\Omega$. The different theorems on elliptic
boundary-value problems were proved by Schechter \cite{Schechter63},
Berezansky, Krein, Roitberg \cite{BerezanskyKreinRoitberg63}, Roitberg
\cite{Roitberg64, Roitberg68, Roitberg96, Roitberg99}, Kostarchuk and Roitberg
\cite{KostarchukRoitberg73} (also see the monograph \cite[Ch.~III,
Sec.~6]{Berezansky68} and the survey \cite[Sec.~7.9]{Agranovich97}). In these
theorems, the solution and/or the right-hand side of the elliptic equation are
not distributions in $\Omega$. One of such theorems proved by Ya.~A.~Roitberg
will be used in this paper (see Proposition~\ref{prop1}).

Ya.~A.~Roitberg \cite[Sec.~2.4]{Roitberg72} considered a condition on the space
$X^{s}(\Omega)$, which was somewhat stronger than our Condition~I$_{s}$. He
demanded additionally that $C^{\infty}(\,\overline{\Omega}\,)\subset
X^{s}(\Omega)$. Under this stronger condition, Roitberg
\cite[Sec.~2.4]{Roitberg72}, \cite[p.~190]{Roitberg96} proved the boundedness
of the operator \eqref{3} for all $s<0$. Remark that Roitberg's condition does
not cover the important case where $X^{s}(\Omega)=\{0\}$ as well as some
weighted spaces $X^{s}(\Omega)=\varrho H^{s}(\Omega)$ which we consider.

We also note that Condition~I$_{s}$ is fulfilled for some classis of the
Hilbertian H\"ormander spaces \cite[Sec.~2.2]{Hermander63},
\cite[Sec.~10.1]{Hermander83}. Their applications to elliptic operators and
elliptic boundary-value problems were studied by V.~A.~Mikhailets and the
author in [22--31] .

The results of this paper are formulated in Section~\ref{sec2} as
Theorems~\ref{th1}, \ref{th2}, \ref{th3}, and
Corollaries~\ref{cor1},~\ref{cor2}. The main result, Theorem~\ref{th1}, is
proved in Section~\ref{sec4}, all the rest in Section~\ref{sec5}. In
Section~\ref{sec3}, we formulate the auxiliary propositions needed for our
proofs. At the end of the paper, we give the appendix, in which a useful
proposition on weight functions of the form $\varrho=\varrho_{1}^{\delta}$ is
established.

\section{Results} \label{sec2}

We introduce some necessary function spaces. Let $s\in\mathbb{R}$. Recall that
$$
H^{s}(\mathbb{R}^{n}):=\bigl\{w\in\mathcal{S}'(\mathbb{R}^{n}):
\,\|w\|_{H^{s}(\mathbb{R}^{n})}:=
\|(1+|\xi|^{2})^{s/2}\,\widehat{w}(\xi)\|_{L_{2}(\mathbb{R}^{n}_{\xi})}<\infty\bigr\}.
$$
Here, $\mathcal{S}'(\mathbb{R}^{n})$ is the topological linear space of
tempered distributions in $\mathbb{R}^{n}$, whereas $\widehat{w}$ is the
Fourier transform of $w$. For a closed set $Q\subset\mathbb{R}^{n}$, we put
$$
H^{s}_{Q}(\mathbb{R}^{n}):=\bigl\{w\in
H^{s}(\mathbb{R}^{n}):\,\mathrm{supp}\,w\subseteq Q\bigr\}.
$$
The space $H^{s}_{Q}(\mathbb{R}^{n})$ is Hilbert with respect to the inner
product in $H^{s}(\mathbb{R}^{n})$. We are interested in the cases where
$Q\in\{\overline{\Omega},\widehat{\Omega},\Gamma\}$ with
$\widehat{\Omega}:=\mathbb{R}^{n}\setminus\Omega$.

Following \cite[Ch.~1, Sec.~12.1]{LionsMagenes68}, we will define the Hilbert
space $H^{s}(\Omega)$. For arbitrary $s\geq0$ we set
$$
H^{s}(\Omega):=H^{s}(\mathbb{R}^{n})/H^{s}_{\widehat{\Omega}}(\mathbb{R}^{n})=
\bigl\{w\!\upharpoonright\!\Omega:\,w\in H^{s}(\mathbb{R}^{n})\bigr\}.
$$
The space $H^{s}(\Omega)$ is complete with respect to the Hilbertian norm
$$
\|u\|_{H^{s}(\Omega)}:=\inf\,\bigl\{\,\|w\|_{H^{s}(\mathbb{R}^{n})}:\, w\in
H^{s}(\mathbb{R}^{n}),\;w=u\;\mbox{in}\;\Omega\,\bigr\}.
$$
The set $C^{\infty}(\,\overline{\Omega}\,)$ is dense in $H^{s}(\Omega)$, each
measurable function $u$ over $\Omega$ being identified with the antilinear
functional $(u,\cdot\,)_{\Omega}$. We denote by $H^{s}_{0}(\Omega)$ the closure
of the linear manifold
$$
C^{\infty}_{0}(\Omega):=\{u\in
C^{\infty}(\,\overline{\Omega}\,):\;\mathrm{supp}\,u\subset\Omega\}
$$
in the topology of $H^{s}(\Omega)$. The space $H^{s}_{0}(\Omega)$ is complete
with respect to the inner product in $H^{s}(\Omega)$.

For arbitrary $s<0$, we denote by $H^{s}(\Omega)$ the Hilbert space antidual to
the space $H^{-s}_{0}(\Omega)$ with respect to the inner product in
$L_{2}(\Omega)$. Since antilinear functionals from $H^{s}(\Omega)$ are defined
uniquely by their values on functions from $C^{\infty}_{0}(\Omega)$, we can
correctly identify these functionals with distributions in $\Omega$. It useful
to keep in mind that
$H^{s}(\Omega)=H^{s}_{\overline{\Omega}}(\mathbb{R}^{n})/H^{s}_{\Gamma}(\mathbb{R}^{n})$
with equality of norms for every $s<0$ \cite[Ch.~1, Remark
12.5]{LionsMagenes68}, and
$H^{s}(\Omega)=H^{s}(\mathbb{R}^{n})/H^{s}_{\widehat{\Omega}}(\mathbb{R}^{n})$
with equivalence of norms for all non half-integer $s<0$ \cite[Sec.
4.8.2]{Triebel78}. It follows from the first equality that
$C^{\infty}_{0}(\Omega)$ is dense in $H^{s}(\Omega)$ for every $s<0$ \cite[Sec.
4.3.2, Theorem~1(b)]{Triebel78}.

Thus, the Hilbert space $H^{s}(\Omega)$ is defined for every $s\in\mathbb{R}$
and imbedded continuously in $\mathcal{D}'(\Omega)$. We have the compact dense
imbedding $H^{s+\delta}(\Omega)\hookrightarrow H^{s}(\Omega)$ if $\delta>0$.

We denote by $H^{s}(\Gamma)$ the Hilbertian Sobolev space with the index
$s\in\mathbb{R}$ defined over the closed compact manifold $\Gamma$ \cite[Ch.~1,
Sec.~7.3]{LionsMagenes68}. The space consists of all distributions on $\Gamma$
belonging to $H^{s}(\mathbb{R}^{n-1})$ in local coordinates and does not depend
on a chose of local charts on $\Gamma$ up to equivalence of norms.

Let us formulate the main result of the paper.

\begin{theorem}\label{th1}
Let $s<0$, and $X^{s}(\Omega)$ be an arbitrary Hilbert space imbedded
continuously in $\mathcal{D}'(\Omega)$ and satisfying Condition
$\mathrm{I}_{s}$. Then the elliptic boundary-value problem \eqref{1} possesses
the following properties:
\begin{itemize}
\item[(i)] The set
$$
D^{\infty}_{A,X}(\Omega):=\{\,u\in C^{\infty}(\,\overline{\Omega}\,):\,Au\in
X^{s}(\Omega)\,\}
$$
is dense in $D^{s+2q}_{A,X}(\Omega)$.
\item[(ii)] The mapping $u\rightarrow(Au,Bu)$ with $u\in D^{\infty}_{A,X}(\Omega)$
can be extended by continuity to the bounded linear operator \eqref{3}.
\item[(iii)] The operator \eqref{3} is Fredholm. Its kernel is $N$, and its range
consists of all vectors
$(f,g_{1},\ldots,g_{q})\in\mathcal{X}_{s}(\Omega,\Gamma)$ satisfying \eqref{4}.
\item[(iv)] If the set $\mathcal{O}(X^{\infty}(\Omega))$ is dense in
$H^{s}_{\overline{\Omega}}(\mathbb{R}^{n})$, then the index of \eqref{3} is
$\dim N-\dim N^{+}$.
\end{itemize}
\end{theorem}

Let us consider some applications of Theorem~\ref{th1} caused by a particular
choice of the space $X^{s}(\Omega)$. Evidently, the space
$X^{s}(\Omega):=\{0\}$ satisfies Condition~I$_{s}$. In this case,
Theorem~\ref{th1} describes properties of the semihomogeneous boundary-value
problem \eqref{1} with $f=0$ and holds true for every $s\in\mathbb{R}$ (also
see \cite{06UMJ11}).

All the Hilbertian Sobolev spaces satisfying Condition I$_{s}$ are found in the
next theorem.

\begin{theorem}\label{th2}
Let $s<0$ and $\sigma\in\mathbb{R}$. The space
$X^{s}(\Omega):=H^{\sigma}(\Omega)$ satisfies Condition~$\mathrm{I}_{s}$ if and
only if
\begin{equation}\label{5}
\sigma\geq\max\{s,-1/2\}.
\end{equation}
\end{theorem}

The next result follows from Theorems \ref{th1} and \ref{th2}.

\begin{corollary}\label{cor1}
Let $s<0$, and \eqref{5} be valid. Then the mapping $u\mapsto(Au,Bu)$ with
$u\in C^{\infty}(\,\overline{\Omega}\,)$ can be extended by continuity to the
bounded and Fredholm operator
\begin{equation}\label{6}
(A,B):\,\bigl\{u\in H^{s+2q}(\Omega):Au\in H^{\sigma}(\Omega)\bigr\}\rightarrow
H^{\sigma}(\Omega)\oplus\bigoplus_{j=1}^{q}\,H^{s+2q-m_{j}-1/2}(\Gamma),
\end{equation}
its domain being the Hilbert space with respect to the norm
$$
\bigl(\,\|u\|_{H^{s+2q}(\Omega)}^{2}+\|Au\|_{H^{\sigma}(\Omega)}^{2}\bigr)^{1/2}.
$$
The index of \eqref{6} is $\dim N-\dim N^{+}$ and independent of $s,\;\sigma$.
\end{corollary}

Here, we note the particular case where $\sigma=s$. If $-1/2<\sigma=s<0$, then
the domain of \eqref{6} coincides with the space $H^{s+2q}(\Omega)$ up to
equivalence of norms. If $\sigma=s=-1/2$, then the domain is narrower than
$H^{2q-1/2}(\Omega)$ but does not depend on $A$ as well (see Sec.~\ref{sec5.3}
below).

We always have $X^{s}(\Omega)\subseteq H^{-1/2}(\Omega)$ in Theorem~\ref{th2}.
But we can get a space $X^{s}(\Omega)$ containing an extensive class of some
distributions $f\notin H^{-1/2}(\Omega)$ and satisfying Condition~I$_{s}$ if we
use certain weighted spaces $\varrho H^{s}(\Omega)$.

A function $\varrho$ given in $\Omega$ is called a \emph{multiplier} in the
space $H^{s}(\Omega)$ if the operator of multiplication by $\varrho$ maps this
space into itself and is bounded on it. A description of the class of all
multipliers in $H^{s}(\Omega)$ with $s\geq0$ was given in
\cite[Sec.~6.3.3]{MazyaShaposhnikova86}.

Let $s<-1/2$. We introduce the following condition on the weight function
$\varrho$.

\medskip

\noindent\textbf{Condition II$_{s}$.} The function $\varrho$ is a multiplier in
the space $H^{-s}(\Omega)$, and
\begin{equation}\label{e}
D_{\nu}^{j}\,\varrho=0\quad\mbox{on}\quad\Gamma\quad\mbox{for every}\quad
j\in\mathbb{Z},\quad0\leq j<-s-1/2.
\end{equation}

\medskip

Here, $D_{\nu}$ is the derivation operator with respect to the unit vector
$\nu$ of inner normal to the boundary $\Gamma$ of $\Omega$. Note that if
$\varrho$ is a multiplier in $H^{-s}(\Omega)$, then evidently $\varrho\in
H^{-s}(\Omega)$. By the trace theorem \cite[Ch.~1, Sec.~9.2]{LionsMagenes68},
there is a trace $(D_{\nu}^{j}\varrho)\!\!\upharpoonright\!\Gamma\in
H^{-s-j-1/2}(\Gamma)$ for every integer $j\geq0$ such that $-s-j-1/2>0$. Hence,
Condition~II$_{s}$ is formulated correctly.

\begin{theorem}\label{th3}
Let $s<-1/2$, and a function $\varrho\in C^{\infty}(\Omega)$ be positive. The
space $X^{s}(\Omega):=\varrho H^{s}(\Omega)$ satisfies
Condition~$\mathrm{I}_{s}$ if and only if the function $\varrho$ satisfies
Condition~$\mathrm{II}_{s}$.
\end{theorem}

The next result follows from Theorems \ref{th1} and \ref{th3}.

\begin{corollary}\label{cor2}
Let $s<-1/2$, and a positive function $\varrho\in C^{\infty}(\Omega)$ satisfy
Condition~$\mathrm{II}_{s}$. Then the mapping $u\rightarrow(Au,Bu)$ with $u\in
C^{\infty}(\,\overline{\Omega}\,)$, $Au\in\varrho H^{s}(\Omega)$ can be
extended by continuity to the bounded and Fredholm operator
\begin{equation}\label{7}
(A,B):\,\bigl\{u\in H^{s+2q}(\Omega):Au\in\varrho
H^{s}(\Omega)\bigr\}\rightarrow \varrho
H^{s}(\Omega)\oplus\bigoplus_{j=1}^{q}\,H^{s+2q-m_{j}-1/2}(\Gamma),
\end{equation}
its domain being the Hilbert space with respect to the norm
$$
\bigl(\,\|u\|_{H^{s+2q}(\Omega)}^{2}+\|\varrho^{-1}Au\|_{H^{s}(\Omega)}^{2}\bigr)^{1/2}.
$$
The index of \eqref{7} is $\dim N-\dim N^{+}$ and independent of $s,\;\varrho$.
\end{corollary}

We give an important example of a function $\varrho$ satisfying
Condition~II$_{s}$ for fixed $s<-1/2$ if we set
$\varrho:=\varrho_{1}^{\delta}$, where $\varrho_{1}$ meets \eqref{f}, and the
number $\delta$ is such that $\delta\geq-s-1/2\in\mathbb{Z}$ or
$\delta>-s-1/2\notin\mathbb{Z}$ (we will prove it in the appendix).

Let us compare Theorem~\ref{th1} and its Corollaries~\ref{cor1}, \ref{cor2}
with the Lions-Magenes theorems \cite{LionsMagenes68, LionsMagenes62V,
LionsMagenes63VI, Magenes65} on elliptic boundary-value problems in the spaces
of distributions.

We restrict ourselves to the case of Hilbertian spaces, whereas the
non-Hilbertian Sobolev spaces were considered in \cite{LionsMagenes62V,
LionsMagenes63VI, Magenes65} as well.

A proposition similar to Theorem~\ref{th1} was proved in
\cite[Sec.~6.10]{Magenes65} for non half-integers $s\leq-2q$ and the Dirichlet
problem, the space $X^{s}(\Omega)$ obeying some different conditions depending
on the problem. Our Condition~I$_{s}$ does not depend on it.

Theorem~LM$_{1}$ is a particular case of Corollary~\ref{cor1} where $\sigma=0$,
i.e. $X^{s}(\Omega)=L_{2}(\Omega)$. Note that some spaces $X^{s}(\Omega)$
containing $L_{2}(\Omega)$ are permissible in Theorem~\ref{th2} and
Corollary~\ref{cor1}. The space $X^{s}(\Omega)=H^{-1/2}(\Omega)$ is the most
extensive among them provided that $s\leq-1/2$. If $-1/2<\sigma=s<0$, then
Corollary~\ref{cor1} coincides with Theorem~7.5 of Lions and Magenes
\cite[Ch.~2]{LionsMagenes68} proved under the additional assumption that
$N=N^{+}=\{0\}$.

Theorem~LM$_{2}$ is a special case of Corollary~\ref{cor2} because the function
$\varrho:=\varrho_{1}^{-s}$ satisfies Condition~II$_{s}$.

\section{Auxiliary propositions} \label{sec3}

First we note \cite[Ch.~1, Theorem 11.5]{LionsMagenes68} that for every $s>1/2$
\begin{equation}\label{8}
H^{s}_{0}(\Omega):=\bigl\{u\in H^{s}(\Omega):\,
D_{\nu}^{j}\,u=0\;\;\mbox{on}\;\;\Gamma\;\;\forall\;\;j\in\mathbb{Z},\;0\leq
j<s-1/2\bigr\}.
\end{equation}
In addition \cite[Ch.~1, Theorem 11.1]{LionsMagenes68}
\begin{equation}\label{9}
H^{s}_{0}(\Omega)=H^{s}(\Omega)\quad\mbox{for}\quad 0\leq s\leq1/2.
\end{equation}

Further we will use some results of Ya.~A.~Roitberg \cite{Roitberg64,
Roitberg96} on properties of the problem \eqref{1} in the Hilbert scale
\begin{equation}\label{10}
H^{s,(2q)}(\Omega):=\widetilde{H}^{s,2,(2q)}(\Omega),\quad s\in\mathbb{R},
\end{equation}
introduced by him. We also need some properties of this scale.

Let us give the definition of the scale \eqref{10}. Let $s\in\mathbb{R}$. In
the case were $s\geq0$ we denote by $H^{s,(0)}(\Omega)$ the space
$H^{s}(\Omega)$. In the case were $s<0$ we denote by $H^{s,(0)}(\Omega)$ the
space $H^{s}_{\overline{\Omega}}(\mathbb{R}^{n})$ antidual to $H^{-s}(\Omega)$
with respect to the inner product in $L_{2}(\Omega)$. The space
$H^{s,(0)}(\Omega)$ is Hilbert for every $s\in\mathbb{R}$, with the set
$C^{\infty}(\,\overline{\Omega}\,)$ being dense in it. Here as usual, the
function $f\in C^{\infty}(\,\overline{\Omega}\,)$ is identified with the
functional $(f,\cdot\,)_{\Omega}$.

In view of \eqref{9}
\begin{equation}\label{11}
H^{s,(0)}(\Omega)=H^{s}(\Omega)\quad\mbox{with equality of
norms}\quad\mbox{for}\quad s\geq-1/2.
\end{equation}
If $s<-1/2$, then the spaces $H^{s,(0)}(\Omega)$ and $H^{s}(\Omega)$ are
different.

Let $s\in\mathbb{R}$ and $s\neq j-1/2$ for all $j=1,\ldots,2q$. We denote by
$H^{s,(2q)}(\Omega)$ the completion of the linear system
$C^{\infty}(\,\overline{\Omega}\,)$ with respect to the norm
$$
\|\,u\,\|_{H^{s,(2q)}(\Omega)}:= \Bigl(\,\|\,u\,\|_{H^{s,(0)}(\Omega)}^{2}+
\sum_{j=1}^{2q}\;\|(D_{\nu}^{j-1}u)\!\!\upharpoonright\!\Gamma\|
_{H^{s-j+1/2}(\Gamma)}^{2}\,\Bigr)^{1/2}.
$$
The space $H^{s,(2q)}(\Omega)$ is separable Hilbert.

In the case where $s\in\{j-1/2:j=1,\ldots,2q\}$ we define the separable Hilbert
space $H^{s,(2q)}(\Omega)$ by means of the complex interpolation
$$
H^{s,(2q)}(\Omega):=\bigl[\,H^{s-1/2,(2q)}(\Omega),
H^{s+1/2,(2q)}(\Omega)\,\bigr]_{1/2}.
$$

We note that by the trace theorem \cite[Ch.~1, Sec.~9.2]{LionsMagenes68}
\begin{equation}\label{12}
H^{s,(2q)}(\Omega)=H^{s}(\Omega)\quad\mbox{with equivalence of
norms}\;\;\mbox{for}\;\;s>2q-1/2.
\end{equation}
The spaces $H^{s,(2q)}(\Omega)$ and $H^{s}(\Omega)$ are different if $s\leq
2q-1/2$.

The imbeddings $H^{s_{2},(0)}(\Omega)\hookrightarrow H^{s_{1},(0)}(\Omega)$ and
$H^{s_{2},(2q)}(\Omega)\hookrightarrow H^{s_{1},(2q)}(\Omega)$ are compact and
dense for arbitrary $s_{1},s_{2}\in\mathbb{R}$, $s_{1}<s_{2}$. This follows
from the compactness of the imbeddings $H^{s_{2}}(\Omega)\hookrightarrow
H^{s_{1}}(\Omega)$ and $H^{s_{2}}(\Gamma)\hookrightarrow H^{s_{1}}(\Gamma)$.

\begin{proposition}[\cite{Roitberg96}, Theorems 4.1.1, 5.3.1] \label{prop1}
Let $s\in\mathbb{R}$. The mapping $u\mapsto(Au,Bu)$ with $u\in
C^{\infty}(\,\overline{\Omega}\,)$ can be extended by continuity to the linear
bounded operator
\begin{equation}\label{13}
(A,B):\,H^{s+2q,(2q)}(\Omega)\rightarrow
H^{s,(0)}(\Omega)\oplus\bigoplus_{j=1}^{q}\,
H^{s+2q-m_{j}-1/2}(\Gamma)=:\mathcal{H}_{s,(0)}(\Omega,\Gamma).
\end{equation}
This operator is Fredholm. Its kernel coincides with $N$, whereas its range
consists of all vectors
$(f,g_{1},\ldots,g_{q})\in\mathcal{H}_{s,(0)}(\Omega,\Gamma)$ satisfying
condition \eqref{4}. The index of \eqref{13} is $\dim N-\dim N^{+}$.
\end{proposition}

Proposition~\ref{prop1} is another example of the generic theorem on elliptic
boundary-value problems. If $s\geq0$, then Proposition~\ref{prop1} coincides
with Theorem~0 by \eqref{12}. But if $s<-1/2$, then both the spaces
$H^{s+2q,(2q)}(\Omega)$ and
$H^{s,(0)}(\Omega)=H^{s}_{\overline{\Omega}}(\mathbb{R}^{n})$ consist of the
elements which are not distributions in the domain $\Omega$.

\begin{proposition}[\cite{Roitberg96}, Theorem 7.1.1] \label{prop2}
Let $s\in\mathbb{R}$, $\delta>0$, and $u\in H^{s+2q,(2q)}(\Omega)$. If
$(A,B)u\in\mathcal{H}_{s+\delta,(0)}(\Omega,\Gamma)$, then $u\in
H^{s+2q+\delta,(2q)}(\Omega)$.
\end{proposition}

\begin{proposition}[\cite{Roitberg96}, Theorem 6.1.1] \label{prop3}
Let $s\in\mathbb{R}$. The following assertions are true:
\begin{itemize}
\item [(i)] The norm in the space $H^{s+2q,(2q)}(\Omega)$ is equivalent to the
norm
\begin{equation}\label{14}
\bigl(\,\|u\|_{H^{s+2q,(0)}(\Omega)}^{2}+
\|Au\|_{H^{s,(0)}(\Omega)}^{2}\bigr)^{1/2}
\end{equation}
on the set of all functions $u\in C^{\infty}(\,\overline{\Omega}\,)$. Therefore
the space $H^{s+2q,(2q)}(\Omega)$ coincides with the completion of the linear
system $C^{\infty}(\,\overline{\Omega}\,)$ with respect to the norm \eqref{14}.
\item [(ii)] The mapping $I_{A}:\,u\mapsto(u,Au)$ with
$u\in C^{\infty}(\,\overline{\Omega}\,)$ can be extended by continuity to the
homeomorphism
$$
I_{A}:\,H^{s+2q,(2q)}(\Omega)\leftrightarrow K_{s+2q,A}(\Omega).
$$
Here,
\begin{align}\notag
K_{s+2q,A}(\Omega):=\bigl\{\, & (u_{0},f):\,u_{0}\in H^{s+2q,(0)}(\Omega),\;
f\in H^{s,(0)}(\Omega), \\ \label{15} &
(u_{0},A^{+}v)_{\Omega}=(f,v)_{\Omega}\;\;\forall\;\;v\in
H^{2q}_{0}(\Omega)\cap H^{-s,(0)}(\Omega)\,\bigr\}
\end{align}
is a closed subspace in $H^{s+2q,(0)}(\Omega)\oplus H^{s,(0)}(\Omega)$.
\end{itemize}
\end{proposition}

\begin{proposition}[\cite{Roitberg96}, Theorem 6.2.1] \label{prop4}
Let $s<-2q-1/2$. For each couple of distributions $u_{0}\in
H^{s+2q,(0)}(\Omega)$ and $f\in H^{s,(0)}(\Omega)$ satisfying the condition
\begin{equation}\label{16}
(u_{0},A^{+}v)_{\Omega}=(f,v)_{\Omega}\quad\mbox{for all}\quad v\in
C^{\infty}_{0}(\Omega),
\end{equation}
there exists a unique couple $(u_{0}^{\ast},f)\in K_{s+2q,A}(\Omega)$ such that
\begin{equation}\label{17}
(u_{0},v)_{\Omega}=(u_{0}^{\ast},v)_{\Omega}\quad\mbox{for all}\quad v\in
C^{\infty}_{0}(\Omega).
\end{equation}
Furthermore,
\begin{equation}\label{18}
\|u_{0}^{\ast}\|_{H^{s+2q,(0)}(\Omega)}\leq c\,
\bigl(\,\|u_{0}\|_{H^{s+2q,(0)}(\Omega)}^{2}+
\|f\|_{H^{s,(0)}(\Omega)}^{2}\bigr)^{1/2},
\end{equation}
with number $c>0$ being independent of $u_{0}$, $f$, and $u_{0}^{\ast}$.
\end{proposition}

\begin{remark}[\cite{Roitberg96}, Sec.~6.2] \label{rem1}
The conditions \eqref{15} and \eqref{16} are equivalent for $s\geq-2q-1/2$, but
they are not equivalent if $s<-2q-1/2$.
\end{remark}

\section{Proof of the main result} \label{sec4}

Now we will prove the main result of the paper, Theorem~\ref{th1}. We assume
that its condition be fulfilled; i.e., $s<0$ and the Hilbert space
$X^{s}(\Omega)$ is imbedded continuously in $\mathcal{D}'(\Omega)$ and
satisfies Condition~I$_{s}$. It follows that the mapping $f\mapsto\mathcal{O}f$
with $f\in X^{\infty}(\Omega)$ can be extended by continuity to the bounded
linear operator
\begin{equation}\label{19}
\mathcal{O}:\,X^{s}(\Omega)\rightarrow
H^{s}_{\overline{\Omega}}(\mathbb{R}^{n})=H^{s,(0)}(\Omega).
\end{equation}
This operator is injective. Indeed, let $\mathcal{O}f=0$ for a distribution
$f\in X^{s}(\Omega)$. Chose a sequence $(f_{k})\subset X^{\infty}(\Omega)$ such
that $f_{k}\rightarrow f$ in
$X^{s}(\Omega)\hookrightarrow\mathcal{D}'(\Omega)$. Then
$\mathcal{O}f_{k}\rightarrow0$ in
$H^{s}_{\overline{\Omega}}(\mathbb{R}^{n})\hookrightarrow\mathcal{S}'(\mathbb{R}^{n})$
that implies
$$
(f,v)_{\Omega}=\lim(f_{k},v)_{\Omega}=\lim(\mathcal{O}f_{k},v)_{\Omega}=0\quad
\mbox{for every}\quad v\in C^{\infty}_{0}(\Omega).
$$
Thus, $f=0$ as a distribution belonging to
$X^{s}(\Omega)\hookrightarrow\mathcal{D}'(\Omega)$; i.~e., the operator
\eqref{19} is injective. This operator defines the continuous imbedding
$X^{s}(\Omega)\hookrightarrow H^{s,(0)}(\Omega)$.

According to Proposition~\ref{prop1}, the element $Au\in H^{s,(0)}(\Omega)$ is
correctly defined for an arbitrary $u\in H^{s+2q,(2q)}(\Omega)$ by means of
passing to the limit. We set
$$
D^{s+2q,(2q)}_{A,X}(\Omega):=\{\,u\in H^{s+2q,(2q)}(\Omega):\,Au\in
X^{s}(\Omega)\,\}.
$$
We also endow the space $D^{s+2q,(2q)}_{A,X}(\Omega)$ with the graphics inner
product
$$
(u_{1},u_{2})_{D^{s+2q,(2q)}_{A,X}(\Omega)}:=(u_{1},u_{2})_{H^{s+2q,(2q)}(\Omega)}+
(Au_{1},Au_{2})_{X^{s}(\Omega)}.
$$
The space $D^{s+2q,(2q)}_{A,X}(\Omega)$ is complete with respect to it. Indeed,
let $(u_{k})$ be a Cauchy sequence in $D^{s+2q,(2q)}_{A,X}(\Omega)$. Since both
the spaces $H^{s+2q,(2q)}(\Omega)$ and $X^{s}(\Omega)$ are complete, the limits
$u:=\lim u_{k}$ in $H^{s+2q,(2q)}(\Omega)$ and $f:=\lim Au_{k}$ in
$X^{s}(\Omega)$ exist. The first of them implies that $Au=\lim Au_{k}$ in
$H^{s,(0)}(\Omega)$. We have from this in view of the second limit and the
continuity of \eqref{19} that $Au=f\in X^{s}(\Omega)$. Hence, $u\in
D^{s+2q,(2q)}_{A,X}(\Omega)$ and $\lim u_{k}=u$ in
$D^{s+2q,(2q)}_{A,X}(\Omega)$; i.e., the space $D^{s+2q,(2q)}_{A,X}(\Omega)$ is
complete.

By Proposition~\ref{prop1}, the restriction of the operator \eqref{13} to
$D^{s+2q,(2q)}_{A,X}(\Omega)$ gives the bounded operator
\begin{equation}\label{20}
(A,B):\,D^{s+2q,(2q)}_{A,X}(\Omega)\rightarrow\mathcal{X}_{s}(\Omega,\Gamma).
\end{equation}
The kernel of \eqref{20} is $N$, and the range consists of all vectors
$(f,g_{1},\ldots,g_{q})\in\mathcal{X}_{s}(\Omega,\Gamma)$ satisfying condition
\eqref{4}. Hence, the operator \eqref{20} is Fredholm, with its co-kernel being
of a dimension $\beta\leq\dim N^{+}$.

Moreover, if $\mathcal{O}(X^{\infty}(\Omega))$ is dense in
$H^{s}_{\overline{\Omega}}(\mathbb{R}^{n})$, then $\beta=\dim N^{+}$. Indeed,
denoting the operator \eqref{13} by $\Lambda$ and the narrower operator
\eqref{20} by $\Lambda_{0}$, we consider the operators $\Lambda^{\ast}$ and
$\Lambda_{0}^{\ast}$ adjoint to them. Since the imbedding
$\mathcal{X}_{s}(\Omega,\Gamma)\hookrightarrow
\mathcal{H}_{s,(0)}(\Omega,\Gamma)$ is continuous and dense, we have
$\ker\Lambda_{0}^{\ast}\supseteq\ker\Lambda^{\ast}$. Hence
$$
\beta=\dim\mathrm{coker}\,\Lambda_{0}=\dim\ker\Lambda_{0}^{\ast}\geq
\dim\ker\Lambda^{\ast}=\dim\mathrm{coker}\,\Lambda=\dim N^{+}.
$$
Therefore $\beta=\dim N^{+}$ and the index of \eqref{20} is equal to $\dim
N-\dim N^{+}$ if $\mathcal{O}(X^{\infty}(\Omega))$ is dense in
$H^{s}_{\overline{\Omega}}(\mathbb{R}^{n})$.

Let us show that the set $D^{\infty}_{A,X}(\Omega)$ is dense in
$D^{s+2q,(2q)}_{A,X}(\Omega)$. Since
$X^{\infty}(\Omega)\times(C^{\infty}(\Gamma))^{q}$ is dense in
$\mathcal{X}_{s}(\Omega,\Gamma)$, we can write by the Gohberg-Krein lemma
\cite[Lemma 2.1]{GohbergKrein60}
\begin{equation}\label{21}
\mathcal{X}_{s}(\Omega,\Gamma)=(A,B)\bigl(D^{s+2q,(2q)}_{A,X}(\Omega)\bigr)\dotplus
\mathcal{Q}(\Omega,\Gamma),
\end{equation}
where $\mathcal{Q}(\Omega,\Gamma)$ is a finite-dimensional subspace satisfying
the condition
\begin{equation}\label{22}
\mathcal{Q}(\Omega,\Gamma)\subset
X^{\infty}(\Omega)\times(C^{\infty}(\Gamma))^{q}.
\end{equation}
Denote by $\Pi$ the projector of the space $\mathcal{X}_{s}(\Omega,\Gamma)$
onto the first term in \eqref{21} parallel to the second term.

Let $u\in D^{s+2q,(2q)}_{A,X}(\Omega)$. Approximate $F:=(A,B)u$ by a sequence
$(F_{k})\subset X^{\infty}(\Omega)\times(C^{\infty}(\Gamma))^{q}$ in the the
topology of $\mathcal{X}_{s}(\Omega,\Gamma)$. We have
\begin{equation}\label{23}
\lim \Pi F_{k}=\Pi F=F=(A,B)u\quad
\mbox{in}\quad\mathcal{X}_{s}(\Omega,\Gamma),
\end{equation}
and by \eqref{22}
\begin{equation}\label{24}
(\Pi F_{k})\subset X^{\infty}(\Omega)\times(C^{\infty}(\Gamma))^{q}
\end{equation}

The Fredholm operator \eqref{20} naturally generates the topological
isomorphism
$$
\Lambda_{0}:=(A,B):\,D^{s+2q,(2q)}_{A,X}(\Omega)/N\leftrightarrow
\Pi(\mathcal{X}_{s}(\Omega,\Gamma)).
$$
In view of \eqref{23},
$$
\lim\Lambda_{0}^{-1}\,\Pi F_{k}=\{u+w:w\in N\}\quad\mbox{in}\quad
D^{s+2q,(2q)}_{A,X}(\Omega)/N.
$$
Hence, there is a seguence of representatives $u_{k}\in
D^{s+2q,(2q)}_{A,X}(\Omega)$ of cosets $\Lambda_{0}^{-1}\,\Pi F_{k}$ such that
\begin{equation}\label{25}
\lim u_{k}=u\quad\mbox{in}\quad D^{s+2q,(2q)}_{A,X}(\Omega).
\end{equation}
In addition, by \eqref{24} we have
$$
(A,B)u_{k}=\Pi F_{k}\in
C^{\infty}(\,\overline{\Omega}\,)\times(C^{\infty}(\Gamma))^{q}.
$$
It follows in view of Proposition~\ref{prop2}, equality \eqref{11} and the
Sobolev imbedding theorem that
$$
u_{k}\in\bigcap_{\delta>0}\,H^{s+2q+\delta,(2q)}(\Omega)=
\bigcap_{\delta>0}\,H^{s+2q+\delta}(\Omega)=C^{\infty}(\,\overline{\Omega}\,).
$$
Thus, in \eqref{25} we have $(u_{k})\subset D^{\infty}_{A,X}(\Omega)$;
therefore $D^{\infty}_{A,X}(\Omega)$ is dense in $D^{s+2q,(2q)}_{A,X}(\Omega)$.

Next, we will consider the cases $-2q-1/2\leq s<0$ and $s<-2q-1/2$ separately.

The first case: $-2q-1/2\leq s<0$. Then $H^{s+2q,(0)}(\Omega)=H^{s+2q}(\Omega)$
in view of \eqref{11}. We use Proposition~\ref{prop3} and consider the mapping
$I_{0}:u\mapsto u_{0}$ in which $u\in D^{s+2q,(2q)}_{A,X}(\Omega)$ and
$(u_{0},f):=I_{A}u$. This mapping establishes the homeomorphism
\begin{equation}\label{26}
I_{0}:\,D^{s+2q,(2q)}_{A,X}(\Omega)\leftrightarrow D^{s+2q}_{A,X}(\Omega).
\end{equation}
Indeed, note that $\eqref{15}\Leftrightarrow\eqref{16}$ for arbitrary $u_{0}\in
H^{s+2q,(0)}(\Omega)=H^{s+2q}(\Omega)$ and $f\in X^{s}(\Omega)\hookrightarrow
H^{s,(0)}(\Omega)$ (see remark~\ref{rem1}). Condition \eqref{16} means that
$Au_{0}=f$ as distributions in~$\Omega$. It follows by Proposition~\ref{prop3}
that $I_{0}(D^{s+2q,(2q)}_{A,X}(\Omega))=D^{s+2q}_{A,X}(\Omega)$. Moreover, we
have the equivalence of the norms:
\begin{align*}
\|u\|_{D^{s+2q,(2q)}_{A,X}(\Omega)}^{2}&=\|u\|_{H^{s+2q,(2q)}(\Omega)}^{2}+
\|f\|_{X^{s}(\Omega)}^{2} \\
&\asymp\|u_{0}\|_{H^{s+2q,(0)}(\Omega)}^{2}+\|f\|_{H^{s,(0)}(\Omega)}^{2}+
\|f\|_{X^{s}(\Omega)}^{2} \\
&\asymp\|u_{0}\|_{H^{s+2q}(\Omega)}^{2}+\|f\|_{X^{s}(\Omega)}^{2}=
\|u_{0}\|_{D^{s+2q}_{A,X}(\Omega)}^{2}.
\end{align*}
Hence, the mapping $I_{0}$ establishes the homeomorphism~\eqref{26}.

It follows from properties of the operator \eqref{20} denoted by $\Lambda_{0}$
and the operator \eqref{26} that
\begin{equation}\label{27}
\Lambda_{0}I_{0}^{-1}:\,D^{s+2q}_{A,X}(\Omega)\rightarrow
\mathcal{X}_{s}(\Omega,\Gamma)
\end{equation}
is a bounded and Fredholm operator, with the range and index being the same as
for \eqref{20}. Since $I_{0}$ bijectively maps the set
$D^{\infty}_{A,X}(\Omega)$ onto itself, this set is dense in
$D^{s+2q}_{A,X}(\Omega)$, and \eqref{27} is an extension by continuity of the
mapping $u\rightarrow(Au,Bu)$ with $u\in D^{\infty}_{A,X}(\Omega)$.
Theorem~\ref{th1} is proved  in the first case.

The second case: $s<-2q-1/2$. Then
$H^{s+2q,(0)}(\Omega)=H^{s+2q}_{\overline{\Omega}}(\mathbb{R}^{n})$. In
addition
\begin{gather}\label{28}
H^{s+2q}(\Omega)=\bigl\{w\!\upharpoonright\!\Omega:\,w\in
H^{s+2q}_{\overline{\Omega}}(\mathbb{R}^{n})\bigr\}, \\ \label{29}
\|u\|_{H^{s+2q}(\Omega)}=\inf\,\bigl\{\,\|w\|_{H^{s+2q}(\mathbb{R}^{n})}:\,
w\in
H^{s+2q}_{\overline{\Omega}}(\mathbb{R}^{n}),\;w=u\;\mbox{in}\;\Omega\,\bigr\}.
\end{gather}
This follows immediately from the equality $H^{s+2q}(\Omega)=
H^{s+2q}_{\overline{\Omega}}(\mathbb{R}^{n})/H^{s+2q}_{\Gamma}(\mathbb{R}^{n})$
mentioned in Section~\ref{sec2}.

Let us denote $Rw:=w\!\upharpoonright\!\Omega$ for
$w\in\mathcal{S}'(\mathbb{R}^{n})$. We will prove that the mapping
$I_{0}:u\mapsto Ru_{0}$ with $u\in D^{s+2q,(2q)}_{A,X}(\Omega)$ and
$(u_{0},f):=I_{A}u$ establishes the topological isomorphism \eqref{26} in the
case under consideration. (In the first case, $Ru_{0}=u_{0}$.) We use
Proposition~\ref{prop3} and note that $\eqref{15}\Rightarrow\eqref{16}$. For an
arbitrary $u\in D^{s+2q,(2q)}_{A,X}(\Omega)$, we have: $Ru_{0}\in
H^{s+2q}(\Omega)$ (see \eqref{28}), $f=Au\in X^{s}(\Omega)$, and
$$
(Ru_{0},A^{+}v)_{\Omega}=(u_{0},A^{+}v)_{\Omega}=(f,v)_{\Omega}\quad\mbox{for
every}\quad v\in C^{\infty}_{0}(\Omega);
$$
i.~e., $ARu_{0}=f$ as distributions in $\Omega$. Therefore $I_{0}u=Ru_{0}\in
D^{s+2q}_{A,X}(\Omega)$. Moreover, in view of \eqref{29} and the definition of
$H^{s+2q,(2q)}(\Omega)$ we have:
\begin{align*}
\|I_{0}u\|_{D^{s+2q}_{A,X}(\Omega)}^{2}&=
\|Ru_{0}\|_{H^{s+2q}(\Omega)}^{2}+\|f\|_{X^{s}(\Omega)}^{2} \\
&\leq \|u_{0}\|_{H^{s+2q}(\mathbb{R}^{n})}^{2}+\|f\|_{X^{s}(\Omega)}^{2}\leq
\|u\|_{D^{s+2q,(2q)}_{A,X}(\Omega)}^{2}.
\end{align*}
Hence, the operator $I_{0}:D^{s+2q,(2q)}_{A,X}(\Omega)\rightarrow
D^{s+2q}_{A,X}(\Omega)$ is bounded.

Now we will show that this operator is bijective. Let $\omega\in
D^{s+2q}_{A,X}(\Omega)$ and $f:=A\omega\in X^{s}(\Omega)$. Due to \eqref{28},
there is a distribution $u_{0}\in H^{s+2q}_{\overline{\Omega}}(\mathbb{R}^{n})$
such that $\omega=Ru_{0}$. Distributions $u_{0}$ and $f$ satisfy condition
\eqref{16} because
$$
(u_{0},A^{+}v)_{\Omega}=(\omega,A^{+}v)_{\Omega}= (f,v)_{\Omega}\quad\mbox{for
every}\quad v\in C^{\infty}_{0}(\Omega).
$$
According to Proposition~\ref{prop4}, for $u_{0}\in H^{s+2q,(0)}(\Omega)$ and
$f\in X^{s}(\Omega)\hookrightarrow H^{s,(0)}(\Omega)$ there exists a unique
couple $(u_{0}^{\ast},f)\in K_{s+2q,A}(\Omega)$ such that condition \eqref{17}
is fulfilled. This implies by Proposition~\ref{prop3} that
$$
u^{\ast}:=I_{A}^{-1}(u_{0}^{\ast},f)\in
D^{s+2q,(2q)}_{A,X}(\Omega),\quad\mbox{and}\quad
I_{0}u^{\ast}=Ru_{0}^{\ast}=Ru_{0}=\omega.
$$
The element $u^{\ast}$ is a unique preimage of $\omega$ in the mapping $I_{0}$.
Indeed, if $I_{0}u'=\omega$ for some $u'\in\nobreak
D^{s+2q,(2q)}_{A,X}(\Omega)$, then the couple $(u'_{0},f'):=I_{A}u'\in
K_{s+2q,A}(\Omega)$ satisfies the following conditions:
$$
f'=ARu'_{0}=A\omega=f\quad\mbox{and}\quad
(u'_{0},v)_{\Omega}=(\omega,v)_{\Omega}=(u_{0},v)_{\Omega}\;\;\forall\;\;v\in
C^{\infty}_{0}(\Omega).
$$
Therefore by Proposition~\ref{prop4}, the couples $(u'_{0},f')=(u'_{0},f)$ and
$(u_{0}^{\ast},f)$ are equal that implies the equality of their preimages $u'$
and $u^{\ast}$ in the mapping $I_{A}$.

Thus, the linear bounded operator \eqref{26} is bijective in the case examined
and therefore is a topological isomorphism by the Banach theorem on inverse
operator. Now using the Fredholm property of \eqref{20} and reasoning as in the
first case, we complete our proof in the second case.

Theorem~\ref{th1} is proved.

\section{Proofs of Theorems \ref{th2}, \ref{th3}, and corollaries} \label{sec5}

\subsection{Proof of Theorem~\ref{th2}}\label{sec5.1} By the condition, $s<0$,
$\sigma\in\mathbb{R}$, and $X^{s}(\Omega):=H^{\sigma}(\Omega)$. Then the set
$X^{\infty}(\Omega)=C^{\infty}(\,\overline{\Omega}\,)$ is dense in the space
$H^{\sigma}(\Omega)$.

\textbf{Sufficiency.}  Let \eqref{5} be fulfilled, i.e.
$\sigma\geq\max\{s,-1/2\}$. Then by \eqref{11} we have
$$
H^{\sigma}(\Omega)=H^{\sigma,(0)}(\Omega)\hookrightarrow H^{s,(0)}(\Omega),
$$
with the imbedding being continuous. We remind that each function $f\in
C^{\infty}(\,\overline{\Omega}\,)$ is iden\-tified with the functional
$(f,\cdot\,)_{\Omega}$, the last being identified with the function
$\mathcal{O}f$ from the space
$H^{s}_{\overline{\Omega}}(\mathbb{R}^{n})=H^{s,(0)}(\Omega)$. Hence
$$
\|\mathcal{O}f\|_{H^{s}(\mathbb{R}^{n})}\leq
c\,\|f\|_{H^{\sigma}(\Omega)}\quad\mbox{for every}\quad f\in
C^{\infty}(\,\overline{\Omega}\,),
$$
with number $c>0$ being independent of $f$. Sufficiency is proved.

\textbf{Necessity.} Let $X^{s}(\Omega):=H^{\sigma}(\Omega)$ satisfy
Condition~I$_{s}$. We assume that $\sigma<0$. (If $\sigma\geq0$, then \eqref{5}
holds true). The operator \eqref{19} establishes the continuous dense imbedding
$H^{\sigma}(\Omega)\hookrightarrow H^{s,(0)}(\Omega)$. This implies that
$$
H^{-s}(\Omega)=(H^{s,(0)}(\Omega))'\subseteq(H^{\sigma}(\Omega))'=
H^{-\sigma}_{0}(\Omega),
$$
we denoting by $H'$ the space antidual to $H$ with respect to the inner product
in $L_{2}(\Omega)$. Hence, $-s\geq-\sigma$. Moreover $-\sigma\leq1/2$, because
if $-\sigma>1/2$, then the function $f\equiv1\in H^{-s}(\Omega)$ would not
belong to $H^{-\sigma}_{0}(\Omega)$ by virtue of \eqref{8}. Thus, $\sigma$
satisfies \eqref{5}. Necessity is proved.

\subsection{Proof of Corollary~\ref{cor1}}\label{sec5.2} Let numbers $s<0$
and $\sigma$ satisfy inequality~\eqref{5}. The boundedness and the Fredholm
property of \eqref{6} follow immediately from Theorems \ref{th1} and \ref{th2}
in which $X^{s}(\Omega):=H^{\sigma}(\Omega)$ and
$D^{\infty}_{A,X}(\Omega)=C^{\infty}(\,\overline{\Omega}\,)$. Moreover, since
the set $\mathcal{O}(C^{\infty}(\,\overline{\Omega}\,))$ identified with
$C^{\infty}(\,\overline{\Omega}\,)$ is dense in
$H^{s}_{\overline{\Omega}}(\mathbb{R}^{n})=H^{s,(0)}(\Omega)$, the index of
\eqref{6} is equal to $\dim N-\dim N^{+}$ by Theorem \ref{th1}(iv) and
therefore independent of $s,\;f$.

\subsection{Remark to Corollary~\ref{cor1}.}\label{sec5.3}
Here we consider the special case of Corollary~\ref{cor1} where $-1/2\leq
\sigma=s<0$. In this case the domain of the operator \eqref{6} does not depend
on $A$. Indeed, if $-1/2<\sigma=s<0$, then we have the bounded operator
$A:H^{s+2q}(\Omega)\rightarrow H^{s}(\Omega)$ because $s$ is not half-integer
\cite[Ch.~1, Proposition 12.1]{LionsMagenes68}. It follows that the domain of
\eqref{6} coincides with $H^{s+2q}(\Omega)$ and therefore does not depend on
$A$.

If $s=\sigma=-1/2$, we cannot reason as above because the space
$H^{-1/2}(\Omega)$ is narrower than $A(H^{2q-1/2}(\Omega))$. However, in view
of Theorem~\ref{th1}(i), equality \eqref{11}, and Proposition~\ref{prop3} we
draw a conclusion that the domain of \eqref{6} is the completion of the set of
all $u\in C^{\infty}(\,\overline{\Omega}\,)$ with respect to the norm
$$
\|u\|_{H^{2q-1/2}(\Omega)}^{2}+\|Au\|_{H^{-1/2}(\Omega)}^{2}=
\|u\|_{H^{2q-1/2,(0)}(\Omega)}^{2}+\|Au\|_{H^{-1/2,(0)}(\Omega)}^{2}\asymp
\|u\|_{H^{2q-1/2,(2q)}(\Omega)}.
$$
Hence, the domain coincides with the space $H^{2q-1/2,(2q)}(\Omega)$
independent of $A$.

\subsection{Proof of Theorem \ref{th3}}\label{sec5.4} By the condition, $s<-1/2$
whereas $\varrho\in C^{\infty}(\Omega)$ is positive. Let us denote by
$M_{\varrho}$ and $M_{\varrho^{-1}}$ the operators of multiplication by
$\varrho$ and $\varrho^{-1}$ respectively. We have the isometric isomomorhism
$M_{\varrho}:H^{s}(\Omega)\leftrightarrow\varrho H^{s}(\Omega)$. It follows
from this and from the density of $C^{\infty}_{0}(\Omega)$ in $H^{s}(\Omega)$
that the set $C^{\infty}_{0}(\Omega)$ is dense in $X^{s}(\Omega):=\varrho
H^{s}(\Omega)$. Hence, the more extensive set $X^{\infty}(\Omega)$ is dense in
$X^{s}(\Omega)$.

We need the following lemma.

\begin{lemma}\label{lem1} Let $s<-1/2$. The multiplication by
$\varrho\in C^{\infty}(\Omega)$ is a bounded operator
\begin{equation}\label{30}
M_{\varrho}:H^{-s}(\Omega)\rightarrow H^{-s}_{0}(\Omega)
\end{equation}
if and only if $\varrho$ satisfies Condition~$\mathrm{II}_{s}$.
\end{lemma}

\begin{proof} \textbf{Necessity.} If the multiplication by $\varrho$ defines
the bounded operator \eqref{30}, then $\varrho$ is a multiplier in
$H^{-s}(\Omega)$ and belongs to $H^{-s}_{0}(\Omega)$. Therefore $\varrho$
satisfies Condition~II$_{s}$ in view of \eqref{8}.

\textbf{Sufficiency.} Let $\varrho$ satisfy Condition~II$_{s}$. We only need to
prove that $\varrho u\in H^{-s}_{0}(\Omega)$ for every $u\in H^{-s}(\Omega)$.
Condition~II$_{s}$ implies that $\varrho\in H^{-s}_{0}(\Omega)$ in view of
\eqref{8}. We chose sequences $(u_{k})\subset
C^{\infty}(\,\overline{\Omega}\,)$ and $(\varrho_{j})\subset
C^{\infty}_{0}(\Omega)$ such that $u_{k}\rightarrow u$ and
$\varrho_{j}\rightarrow\varrho$ в $H^{-s}(\Omega)$. Since both the functions
$\varrho$ and $u_{k}$ are multipliers in the space $H^{-s}(\Omega)$, we have
therein
$$
\lim_{k\rightarrow\infty}(\varrho u_{k})=\varrho u\quad\mbox{and}\quad
\lim_{j\rightarrow\infty}(\varrho_{j}u_{k})=\varrho u_{k}\quad\mbox{for every
$k$}.
$$
This in view of $\varrho_{j}u_{k}\in C^{\infty}_{0}(\Omega)$ implies that
$\varrho u\in H^{-s}_{0}(\Omega)$. Sufficiency is proved.
\end{proof}

Now let us define the following space with inner product:
$$
\varrho^{-1}H^{-s}_{0}(\Omega):=\{f=\varrho^{-1}v:\,v\in
H^{-s}_{0}(\Omega)\,\},
$$
$$
(f_{1},f_{2})_{\varrho^{-1}H^{-s}_{0}(\Omega)}:= (\varrho f_{1},\varrho
f_{2})_{H^{-s}(\Omega)}.
$$
We have the isometric isomorphism
\begin{equation}\label{31}
M_{{\varrho}^{-1}}:H^{-s}_{0}(\Omega)\leftrightarrow\varrho^{-1}H^{-s}_{0}(\Omega).
\end{equation}
Hence, the space $\varrho^{-1}H^{-s}_{0}(\Omega)$ is complete, with
$C^{\infty}_{0}(\Omega)$ being dense in it.

Note that
\begin{equation}\label{32}
(\varrho^{-1}H^{-s}_{0}(\Omega))'=\varrho H^{s}(\Omega)\quad\mbox{with equality
of norms}.
\end{equation}
Indeed, passing in \eqref{31} to adjoint operator, we get the isometric
isomorphism
$$
M_{{\varrho}^{-1}}:(\varrho^{-1}H^{-s}_{0}(\Omega))'\leftrightarrow
(H^{-s}_{0}(\Omega))'=H^{s}(\Omega).
$$
This by the definition of $\varrho H^{s}(\Omega)$ implies the isometric
isomorphism
$$
I=M_{\varrho}M_{{\varrho}^{-1}}:(\varrho^{-1}H^{-s}_{0}(\Omega))'\leftrightarrow\varrho
H^{s}(\Omega),
$$
where $I$ is the identity operator. Thus, \eqref{32} is proved.

Now we can complete the proof of Theorem~\ref{th3} in the following way.
According to Lemma~\ref{lem1}, Condition~II$_{s}$ is equivalent to the
boundedness of the operator \eqref{30} that by \eqref{31} is equivalent to the
continuous imbedding $H^{-s}(\Omega)\hookrightarrow
\varrho^{-1}H^{-s}_{0}(\Omega)$. This imbedding is dense and by \eqref{32} is
equivalent to the continuous dense imbedding
$$
\varrho H^{s}(\Omega)=(\varrho^{-1}H^{-s}_{0}(\Omega))'
\hookrightarrow(H^{-s}(\Omega))'=H^{s}_{\overline{\Omega}}(\mathbb{R}^{n}).
$$
Finally, the imbedding $\varrho H^{s}(\Omega)\hookrightarrow
H^{s}_{\overline{\Omega}}(\mathbb{R}^{n})$ is equivalent to Condition~I$_{s}$.
Note that the last imbedding is dense because $C^{\infty}_{0}(\Omega)$ is dense
in $H^{s}_{\overline{\Omega}}(\mathbb{R}^{n})$ \cite[Sec.~4.3.2,
Theorem~1(b)]{Triebel78}. Thus, Conditions II$_{s}$ and I$_{s}$ are equivalent
for $X^{s}(\Omega)=\varrho H^{s}(\Omega)$.

Theorem~\ref{th3} is proved.

\subsection{Proof of Corollary~\ref{cor2}}\label{sec5.5}
Let $s<-1/2$, and a positive function $\varrho\in C^{\infty}(\Omega)$ satisfy
Condition~$\mathrm{II}_{s}$. The boundedness and Fredholm property of the
operator \eqref{7} follow from Theorems \ref{th1} and \ref{th3} with
$X^{s}(\Omega):=\varrho H^{s}(\Omega)$. In addition, since the set
$\mathcal{O}(X^{s}(\Omega))$ is dense in
$H^{s}_{\overline{\Omega}}(\mathbb{R}^{n})$, Theorem ~\ref{th1}(iv) implies
that the index of \eqref{7} is $\dim N-\dim N^{+}$ and independent of $s$,
$\varrho$.

\medskip

\noindent\textit{The author would like to thank Yu.~M.~Berezansky and
V.~A.~Mikhailets for their valuable remarks and interest in this work.}

\section*{Appendix}

\subsection*{\textmd{A.1}}\label{sec6.1}

In the appendix we will prove the following proposition, which gives an
important example of the function satisfying Condition~II$_{s}$.

\medskip

\noindent\textbf{Proposition A.} \it Let a number $s<-1/2$ and a function
$\varrho_{1}$ satisfying condition~\eqref{f} be given. Assume that
$\delta\geq-s-1/2\in\mathbb{Z}$ or $\delta>-s-1/2\notin\mathbb{Z}$. Then the
function $\varrho:=\varrho_{1}^{\delta}$ satisfies Condition~$\mathrm{II}_{s}$.
\rm

\subsection*{\textmd{A.2.} Proof of Proposition A}

Condition~\eqref{e} is fulfilled for the function
$\varrho=\varrho_{1}^{\delta}$ because $\varrho_{1}=0$ on $\Gamma$, and
$\delta\geq-s-1/2$. Therefore we only need to prove that $\varrho_{1}^{\delta}$
is a multiplier in the space $H^{-s}(\Omega)$. If the positive number $\delta$
is integer, then the function $\varrho_{1}^{\delta}$ belongs to
$C^{\infty}(\overline{\Omega})$ and therefore is a multiplier in
$H^{-s}(\Omega)$. Further we assume that $\delta\notin\mathbb{Z}$. Then by the
condition, $\delta>-s-1/2$.

It is not difficult to verify that the function $\eta_{\delta}(t):=\nobreak
t^{\delta}$, $0<t<1$, belongs to $H^{-s}((0,1))$ (we will do it in the next
subsection). Hence, this function has an extension from the interval $(0,1)$ to
$\mathbb{R}$ pertaining to $H^{-s}(\mathbb{R})$. Let us retain the notation
$\eta_{\delta}$ for the extension. By the Strichartz theorem
\cite{Strichartz67}, \cite[Sec.~2.2.9]{MazyaShaposhnikova86}, every function
from the space $H^{-s}(\mathbb{R})$ is a multiplier therein if $-s>1/2$. Hence,
$\eta_{\delta}$ is a multiplier in $H^{-s}(\mathbb{R})$. Then \cite[Sec.~2.4,
Proposition~5]{MazyaShaposhnikova86} the function
$\eta_{\delta,n}(t',t_{n}):=\eta_{\delta}(t_{n})$ of arguments
$t'\in\mathbb{R}^{n-1}$, $t_{n}\in\mathbb{R}$ is a multiplier in
$H^{-s}(\mathbb{R}^{n})$. This function coincides with $\varrho_{1}^{\delta}$
in the special local coordinates $(x',t_{n})$ near the boundary $\Gamma$. Here,
$x'$ is a coordinate of a point on $\Gamma$, and $t_{n}$ is the distance from
$\Gamma$. It follows by \cite[Sec.~6.4.1, Lemma~3]{MazyaShaposhnikova86} that
$\varrho_{1}^{\delta}$ is a multiplier in each space $H^{-s}(\Omega\cap
V_{j})$, where $\{V_{j}:j=1,\ldots,r\}$ is a finite collection of balls in
$\mathbb{R}^{n}$ with a sufficiently small radius $\varepsilon$, and the
collection covers the boundary $\Gamma$. By supplementing this collection with
the set $V_{0}:=\{x\in\Omega:\mathrm{dist}(x,\Gamma)>\varepsilon/2\}$, we get
the finite open covering of the closed domain $\overline{\Omega}$. Let certain
functions $\chi_{j}\in C^{\infty}_{0}(V_{j})$, $j=0,1,\ldots,r$, form the
partition of unity on $\overline{\Omega}$ corresponding to this covering. Since
the multiplication by a function from $C^{\infty}_{0}(V_{j})$ is a bounded
operator in the space $H^{-s}(\Omega\cap V_{j})$, the function
$\chi_{j}\varrho_{1}^{\delta}$ has to be a multiplier in this space and
therefore in $H^{-s}(\Omega)$. Hence,
$\varrho_{1}^{\delta}=\sum_{j=0}^{r}\chi_{j}\varrho_{1}^{\delta}$ is a
multiplier in $H^{-s}(\Omega)$.

It remains to proof that $\eta_{\delta}\in H^{-s}((0,1))$. We use the inner
description of the space $H^{-s}((0,1))$. If $-s\in\mathbb{Z}$, the inclusion
$\eta_{\delta}\in H^{-s}((0,1))$ is equivalent to that $\eta_{\delta}\in
L_{2}((0,1))$ and $\eta_{\delta}^{(-s)}\in L_{2}((0,1))$. The last two
inclusions are fulfilled because $\delta>-s-1/2$. Hence, $\eta_{\delta}\in
H^{-s}((0,1))$ in the case examined. If $-s\notin\mathbb{Z}$, then by
\cite[p.~214, Sec.~7.48]{Adams75} the inclusion $\eta_{\delta}\in
H^{-s}((0,1))$ is equivalent to that $\eta_{\delta}\in H^{[-s]}((0,1))$ and
\begin{equation}\label{6.1}
\int_{0}^{1}\int_{0}^{1}\frac{|D_{t}^{[-s]}t^{\delta}-D_{\tau}^{[-s]}\tau^{\delta}|^{2}}
{|t-\tau|^{1+2\{-s\}}}\,dt\,d\tau<\infty.
\end{equation}
Here as usual, $[-s]$ and $\{-s\}$ are the integral and fractional parts of
$-s$ respectively. Since $\delta>[-s]-1/2$, we have $\eta_{\delta}\in
H^{[-s]}((0,1))$, that was proved above. In addition, inequality \eqref{6.1}
holds true by virtue of the following elementary lemma, which we will prove in
Subsection~A.3.

\medskip

\noindent\textbf{Lemma A.} \it Let $\alpha,\beta,\gamma\in\mathbb{R}$, and in
addition, $\alpha\neq0$, $\gamma>0$. Then
\begin{equation}\label{6.2}
I(\alpha,\beta,\gamma):=\int_{0}^{1}\int_{0}^{1}\frac{|t^{\alpha}-\tau^{\alpha}|^{\gamma}}
{|t-\tau|^{\beta}}\,dt\,d\tau<\infty
\end{equation}
if and only if the following inequalities are fulfilled:
\begin{equation}\label{6.3}
\alpha\gamma-\beta>-2,\quad\gamma-\beta>-1,\quad\alpha\gamma>-1.
\end{equation} \rm

\medskip

Indeed, the double integral in \eqref{6.1} is equal to
$c\,I(\alpha,\beta,\gamma)$, where $c$ is a positive number, whereas
$\alpha=\delta-[-s]$, $\beta=1+2\{-s\}$, and $\gamma=2$. Equalities \eqref{6.3}
are fulfilled for these numbers $\alpha$, $\beta$, and $\gamma$, because
$$
\alpha\gamma-\beta=2(\delta+s)-1>-2,\quad\gamma-\beta=1-2\{-s\}>-1,\quad
\alpha\gamma=2(\delta-[-s])>-1.
$$
We have used the condition $\delta>-s-1/2$ in the first and the third
inequalities. Thus, the inclusion $\eta_{\delta}\in H^{-s}((0,1))$ is also
valid in the case of non-integer $s<-1/2$.

Proposition A is proved.

\subsection*{\textmd{A.3.} Proof of Lemma~A}

Changing the variable $\lambda:=\tau/t$ in the inner integral, we can write the
following in view of evident transformations:
$$
I(\alpha,\beta,\gamma)=2\int_{0}^{1}dt\int_{0}^{t}\frac{|t^{\alpha}-\tau^{\alpha}|^{\gamma}}
{|t-\tau|^{\beta}}\,d\tau=
2\int_{0}^{1}t^{\alpha\gamma-\beta+1}dt\int_{0}^{1}\frac{|1-\lambda^{\alpha}|^{\gamma}}
{|1-\lambda|^{\beta}}\,d\lambda.
$$
Here, the integral in variable $t$ is finite if and only if
$\alpha\gamma-\beta>-2$, whereas the integral in $\tau$ is finite if and only
if both the inequalities $\alpha\gamma>-1$ and $\gamma-\beta>-1$ hold. Hence
$\eqref{6.2}\Leftrightarrow\eqref{6.3}$, which is what had to prove.

\end{document}